\documentclass[12pt,a4paper]{amsart}

\usepackage {amsfonts}
\usepackage[english]{babel}
\def\bd{{\mathbb{D}}}

\def\br{{\mathbb{R}}}

\def\bt{{\mathbb{T}}}

\def\ca{{\mathcal{A}}}
\def\css{{\mathcal{S}}}
\def\cc{{\mathcal{C}}}
\def\ch{{\mathcal{H}}}

\def\a{\alpha}
\def\b{\beta}
\def\l{\lambda}

\def\s{\sigma}

\def\p{\varphi}

\def\th{\theta}
\def\d{\delta}

\def\o{\omega}
\def\O{\Omega}

\def\R{{\mathbb R}}

\def\N{{\mathbb N}}
\def\D{{\mathbb D}}
\def\T{{\mathbb T}}

\def\g{\gamma}
\def\G{\Gamma}
\def\ep{\varepsilon}
\def\z{\zeta}

\def\ovl{\overline}

\def\ti{\tilde}

\def\lp{\left(}
\def\rp{\right)}
\def\dist{{\rm dist}}

\def\bs{~\hfill\rule{7pt}{7pt}}

\newtheorem{theorem}{Theorem}
\newtheorem*{Th}{Theorem BGK}

\begin{document}

\centerline{To V.P.~Havin on the occasion of his 75th birthday}

\vskip2cm

\title[A Blaschke-type condition and application]
{A Blaschke-type condition for analytic and subharmonic functions
and application to contraction operators}
\author[S. Favorov and L. Golinskii]{S. Favorov and L. Golinskii}

\address{Mathematical School, Kharkov National University, Swobody sq.4,
Kharkov, 61077 Ukraine}

 \email{Sergey.Ju.Favorov@univer.kharkov.ua}

\address{Mathematics Division, Institute for Low Temperature Physics and
Engineering, 47 Lenin ave., Kharkov 61103, Ukraine} \email{leonid.golinskii@gmail.com}



\begin{abstract}
Let $E$ be a closed set on the unit circle. We find a Blaschke-type condition, optimal
in a sense of the order, on the Riesz measure of a subharmonic function $v$ in the
unit disk with a certain growth at the direction of $E$. In particular case when $E$
is a finite set, and $v=\log|f|$ with an analytic function $f$, our result agrees with
the recent one by A. Borichev, L. Golinskii and S. Kupin. An application to
contractions close to unitary operators in the Hilbert space is given.

\end{abstract}

\maketitle

\vspace{-0.5cm}
\section{Introduction}
\label{s0}  In \cite{bgk} the authors study zero sets of analytic
functions in the unit disk $\bd$, which grow exponentially fast near
a finite set $E$ of points on the unit circle $\bt$. Here is the
main result of \cite{bgk}. As usual, $x_+=\max\{x,0\}$.

\begin{Th}\label{BGK} Let $E\subset\bt$ be a finite set,
$f\in\ca(\bd)$ an analytic function in $\bd$, $|f(0)|=1$, and

\begin{equation}\label{bgk1}
|f(z)|\le \exp\Bigl(\frac{D}{\rho^q(z,E)}\Bigr) \end{equation} with $D,q\ge 0$,
$\rho(z,E)=\dist (z,E)$. Let $Z_f$ be the zero set of $f$, each zero is counted
according to its multiplicity. Then for any $\ep>0$,
\begin{equation}\label{bgk2}
\sum_{z\in Z_f} (1-|z|)\,\rho(z,E)^{(q-1+\ep)_+}\le C(\ep,q,E)\, D.
\end{equation}
\end{Th}

Note that both \eqref{bgk1} and \eqref{bgk2} make sense for arbitrary (infinite)
subsets $E$ on $\bt$, so the question arises naturally, whether the finiteness of $E$,
which is an essential ingredient in the proof of Theorem BGK, can be relaxed. We
should also mention that in the case $E=\bt$, $\rho(z,E)=1-|z|$ we come to the well
known Blaschke-type condition for analytic functions in $\bd$ with radial growth. The
theory of such functions goes back to M.M. Djrbashian \cite{Dj}, see also V.I.~Matsaev
\& E.Z.~Mogulskii \cite{mamo}, W.~Hayman \& B.~Korenblum \cite{hk}, F.A.~Shamoyan
\cite{Sh2}, A.M.~Jerbashian \cite{J}.

Our investigation shows that the natural setting for the above
problem is the class of subharmonic functions $v$ and their Riesz
measures (generalized Laplacians) $\mu=(1/2\pi)\triangle v$ rather
than analytic functions and their zero sets. It turns out that
arbitrary closed sets $E$ can be involved, and the results are
optimal. Note that in the case of subharmonic functions of the form
$\log|f|$ with $f\in\ca(\bd)$, the Riesz measure is a discrete and
integer-valued measure supported on $Z_f$, and $\mu\{z\}$ equals the
multiplicity of zero of $f$ at $z$.

To formulate our result we need the following quantitative
characteristic of the ``sparseness'' of a closed set $E$ in terms of
its $t$-neighborhood
\begin{equation}\label{neib}
 E_t:=\{\z\in\bt: \rho(\z,E)<t\}.
\end{equation}
Put
\begin{equation}\label{aset}
I(\a,E):=\int_0^2\,\frac{\nu_E(s)}{s^{\a+1}}\,ds\le +\infty, \qquad
\nu_E(t)=|E_t|, \quad \a\in\br.
\end{equation}
Here and in what follows we denote by $|A|$ the normalized Lebesgue
measure of a set $A\subset\bt$. It is not hard to see that
$I(\a,E)<\infty$ for any $E$ and an arbitrary negative $\a$. In the
case when $\nu_E(t)=O(t^\b)$, $\b>0$, one has $I(\a,E)<\infty$ for
$\a<\b$. Examples of sets $E$ and evaluations of their
characteristic $I$ are given in Section \ref{s1}.

The main result of our paper is the following

\begin{theorem}\label{t1}  Let $E\subset\bt$ be a closed set.
Let $v$ be a subharmonic function in $\bd$, $v\not\equiv -\infty$,
$\mu$ its Riesz measure, and for all $z\in\D$ and some $q>0$
 $$
v(z)\le \frac{1}{\rho^q(z,E)}\,.
 $$
$(i)$ If $I(q,E)<\infty$, then
\begin{equation}\label{blas1}
\int_\bd (1-|\l|)\,d\mu(\l)<\infty.
\end{equation}
Moreover, if $v(0)\ge 0$, then
\begin{equation}\label{blas2}
\int_\bd (1-|\l|)\,d\mu(\l)\le 1+q2^q I(q,E).
\end{equation}
$(ii)$ If $I(q,E)=\infty$ and $I(\a,E)<\infty$ for some $\a<q$, then
\begin{equation}\label{blas3}
\int_\bd (1-|\l|)\lp\rho(\l,E)\rp^{q-\a}d\mu(\l)<\infty.
\end{equation}
Moreover, if $v(0)\ge 0$, then
\begin{equation}\label{blas4}
\int_\bd (1-|\l|)\lp\rho(\l,E)\rp^{q-\a}d\mu(\l)<36(2^{-\a}+q(90)^{q-\a}\,I(\a,E)).
\end{equation}
\end{theorem}

The result turns out to be optimal in the following sense.

\begin{theorem}\label{t1.1} Let $E\subset\bt$ be a closed set so that
$I(\a,E)=+\infty$ for some $\a\ge0$. Then for the subharmonic
function $v_0(z)=\rho^{-q}(z,E)$ with the Riesz measure $\mu_0$, and
$q\ge\a$, we have
 \begin{equation}\label{c}
  \int_\bd
(1-|\l|)\lp\rho(\l,E)\rp^{q-\a}d\mu_0(\l)=+\infty.
 \end{equation}
\end{theorem}

As a consequence we obtain
\begin{theorem}\label{t2} Let $E\subset\bt$ be a closed set, and
$I(\a,E)<\infty$ for some $\a\in\br$. Let $f\in\ca(\bd)$, and
\begin{equation}\label{e010}
|f(z)|\le \exp\Bigl(\frac{D}{\rho^q(z,E)}\Bigr), \quad z\in\D.
\end{equation} Then
\begin{equation}\label{e011} \sum_{z\in Z_f} (1-|z|)\,\rho(z,E)^{(q-\a)_+}<\infty.
\end{equation}
 Moreover, if $|f(0)|\ge 1$, then
 $$
\sum_{z\in Z_f} (1-|z|)\,\rho(z,E)^{(q-\a)_+}\le C(q,\a,E)\,D,
 $$
$C(q,\a,E)$ does not depend on $f$. In the case $q>\a$ one can take
$C(q,\a,E)=36(2^{-\a}+q(90)^{q-\a}\,I(\a,E))$.
\end{theorem}

The classical Blaschke condition arises for $E=\T$, $\a=-q<0$ and
$q\to 0$.

The first statement of Theorem \ref{t2} is a certain uniqueness
theorem: under assumption \eqref{e010} the divergence of
\eqref{e011} yields $f\equiv 0$.

We proceed as follows. In Section \ref{s1}  the main results are
proved. In Section \ref{s2} an application to the spectral theory of
contractions in the Hilbert space is discussed.

\section{Proof of the main results}\label{s1}

There is a simple way to compute the function $\nu_E$ defined in
\eqref{aset}, in terms of the complimentary arcs of $E$. \footnote{
An interesting case occurs when the number of complimentary arcs is
infinite.} Let
$$
\T\setminus E=\bigcup_j\,\g_j, \qquad|\g_j|\downarrow0.
$$ Then
$$
\nu_E(t)=\sum_{j=N+1}^\infty |\g_j|+\frac{2N}{\pi}\arcsin\frac{t}2+|E|,
 $$
  where $N=N(t)$ is taken from
 $$
|\g_{N+1}|\le\frac2{\pi}\,{\arcsin\frac{t}2}<|\g_N|.
 $$
So the exponent $\a$ in \eqref{aset}, for which $I(\a,E)<\infty$ can
be easily determined. For instance, for
$$
E=\left\{e^{i\p_n}: \p_n=\sum_{k=n}^\infty \frac1{2^k} \right\}\cup \{1\}
$$ one can take any $\a<1$, for
$$
E=\left\{e^{i\p_n}: \p_n=c\sum_{k=n}^\infty \frac1{k^\g}, \quad \g>1\right\}\cup\{1\},
\quad c^{-1}=\sum_{k=1}^\infty \frac1{k^\g}\,,
$$ any $\a<1-\frac1{\g}$ works, for
$$
E=\left\{e^{i\p_n}: \p_n=c\sum_{k=n}^\infty \frac1{k\log^2 k}\right\}\cup\{1\}, \quad
c^{-1}=\sum_{k=2}^\infty \frac1{k\log^2 k}\,,
$$
$I(\a,E)<\infty$ for any $\a<0$. For the generalized Cantor set
$\cc_\b$ (for the standard Cantor set $\b=1/3$) one can take any
$\a<1-d(\b)$,
$$ d(\b)=\frac{\log 2}{\log 2-\log(1-\b)} $$
is the Hausdorff dimension of $\cc_\b$  (see, e.g., \cite{f}).

Let us list some elementary properties of $\nu_E$.

\begin{enumerate}
    \item $\nu_E$ is a continuous and strictly monotone increasing
    function on $[0,t_0(E)]$ and $\nu_E(t)=1$ for $t_0(E)\le t<\infty$;
    \item $\nu_E(0+)=|E|$;
    \item  $\nu_E(t)\ge |E|+t/\pi$ for each closed $E\not=\T$ and $0<t\le
    t_1(E)$;
    \item $\nu_E(t)=O(t)$, as $t\to0$, if and only if $E$ is a finite set.
\end{enumerate}

\medskip

Our argument relies upon some basic results from the potential
theory in the complex plane (see \cite[Chapters 3,4]{ran}). Let $v$
be a subharmonic function on $\bd$, $v\not\equiv -\infty$. The Riesz
measure $\mu=(1/2\pi)\triangle v$ is known to be a positive Radon
measure on $\bd$. If $v$ has a harmonic majorant on a subdomain
$\O\subseteq\bd$, then the following representation holds
\begin{equation}\label{Rie}
  v(z)=u(z)-\int_\O G_\O(z,\l)~d\mu(\l),\quad z\in \O,
\end{equation}
where $u$ is the least harmonic majorant in $\O$, and $G_\O(z,\l)$
is Green's function for $\O$, that is,
 \begin{equation}\label{G}
G_\O(z,\l)=\log(1/|z-\l|)-h(z,\l),
 \end{equation}
$h(z,\cdot)$ is a harmonic function in $\O$ such that
$h(z,\l)=\log(1/|z-\l|)$ for $\l\in\partial\O$. Note that if
$\O\not=\D$ and $0\in\O$, then $h(0,\l)>0$ for $\l\in\O$.

\bigskip
{\it Proof of Theorem \ref{t1}}.  We begin with a proof of
(\ref{blas2}). It easily follows from $I(q,E)<\infty$, $q>0$, and
\eqref{aset} that now $\nu_E(0+)=|E|=0$. Put $\rho(z)=\rho(z,E)$.

The monotone permutation theorem (also known as the ``layer cake
representation'', see, e.g., \cite[Theorem 1.13]{lilo})
\begin{equation}\label{lcr}
\int_X\,f^r(x)\,d\sigma(x)=
r\int_0^\infty\,y^{r-1}\,\sigma(\{x:f(x)>y\})\,dy,
\end{equation}
is of importance in our argument. Here $r>0$, $(X,\sigma)$ is a
measure space, $f\ge 0$ a measurable function on $X$. We have (see
\eqref{aset})
\begin{equation}\label{norm}
 \int_\bt\frac{dm(\z)}{\rho^q(\z)} = q\int_0^\infty
y^{q-1} \left|\left\{\z:\rho(\z)<\frac1y\right\}\right|\,dy
\end{equation}
 $$
 =
2^{-q}+q I(q,E)<\infty,
 $$
$dm$ is the normalized Lebesgue measure on $\bt$. Therefore the
Poisson integral
$$
U(z)=\int_{\bt}\,\frac{1-|z|^2}{|\z-z|^2}\,\frac{dm(\z)}{\rho^q(\z)}\,,
\quad z\in\D
$$
is a well defined harmonic function.

Let us prove that $v$ admits a harmonic majorant on the whole disk
$\bd$. A simple geometrical inequality shows that for all $|z|\le 1$
and $0<\tau<1$ one has $\rho(z)\le 2\rho(\tau z)$. The function
$v_0=\rho^{-q}$ is obviously subharmonic in $\bd$ since $$
v_0(z)=\sup_{\z\in E}(|z-\z|^{-q})$$ (see \cite[Theorem
2.4.7]{ran}), and continuous at any point $\z\in\bt\setminus E$. By
the well known property of the Poisson integral
$$
 \lim_{z\to\z} U(z)=v_0(\z)\ge 2^{-q}v_0(\tau\z),
\qquad \z\in\bt\setminus E, \quad 0<\tau<1.
$$
As long as $\z\in E$, one has $\lim_{z\to\z} U(z)=+\infty$. But
$v_0(\tau z)$ is a {\it bounded} subharmonic function in $\bd$, so
by the Maximum Principle $U(z)\ge 2^{-q} v_0(\tau z)$, $z\in\bd$. It
remains only to tend $\tau\to 1$ to make sure that $2^qU$ is a
desired harmonic majorant for $v$.

Representation \eqref{Rie} takes now the form
 $$
  v(z)=u(z)-\int_{\bd}\log\left|{1-\bar\l z\over
  z-\l}\right|~d\mu(\l),\quad z\in\bd,
 $$
where $u\le U$ in $\bd$. For $z=0$ we come to
 $$
\int_{\bd}(1-|\l|)~d\mu(\l)\le \int_{\bd}\log{1\over |\l|}~d\mu(\l)
\le u(0)\le 2^q U(0)=2^q\|\rho^{-q}\|_{L^1(\T)},
 $$
and \eqref{blas2} follows from \eqref{norm}.

\bigskip

 The case when $I(q,E)=\infty$ but
$I(\a,E)<\infty$ for some $\a<q$ is much more delicate.

For fixed $t\in(0,1)$ denote by $\O=\O_t$ the connected component of
the open subset $\{z\in\D:\,\rho(z)>t\}$, which contains the origin,
so $\O$ is a nonempty subdomain of $\D$. Put
 $$
 \G:=\{z\in\D:\,\rho(z)=t\},\qquad E_t^c=\T\setminus E_t,
  $$
$E_t$ is defined in \eqref{neib}. It is easy to check that $E_t$ is
a finite union of disjoint open arcs, so $E_t^c$ is a finite union
of disjoint closed arcs.  The boundary $\partial\O$ is contained in
the closure of $E_t^c\cup\G$. The function $v$ is bounded from above
in $\O$, so \eqref{Rie} holds for $v$ in $\O$, with the least
harmonic majorant $u$.

Let $V$ be a harmonic function  in $\D$ such that
\begin{equation}\label{harm}
V(\z)=
\left\{%
\begin{array}{ll}
    t^{-q}, & \hbox{$\z\in E_t$,} \\
    \rho^{-q}(\z), & \hbox{$\z\in E_t^c$.} \\
\end{array}%
\right.
\end{equation}
Note that $V$ is continuous in the closed unit disk. For $0<t<1$ and
each $z\in\G$ there is $\z'\in E$ such that $|z-\z'|=t$. Let
$\o(\l,\g_z,\D)$ be the harmonic measure of the arc
$\g_z=\{\z\in\bt: |\z-\z'|\le t\}$. It easily follows from the
explicit expression (see, e.g., \cite[Chapter 1]{ga}) that
$\o(\l,\g_z,\D)\ge1/6$ for $|\l-\z'|=t$, $\l\in\D$. Since
$\g_z\subset E_t$, we get $V(z)\ge t^{-q}\o(z,\g_z,\D)\ge t^{-q}/6$.
Therefore, we have
$$
\limsup_{z\to\z} v(z)\le 6V(\z), \qquad \forall\,\z\in\partial\O.
 $$
By the Maximum Principle
 $$
 v(z)\le 6V(z)\quad\hbox{and so}\quad u(z)\le 6V(z), \qquad z\in\O.
 $$

Let us  prove \eqref{blas4}. Since $v(0)\ge 0$, \eqref{Rie} for
$z=0$ can be written as
\begin{equation}\label{e14}
\int_{\O} G_{\O}(0,\l)\,d\mu(\l)\le u(0)\le
6V(0)=\frac{6\nu_E(t)}{t^q}+6\int_{E_t^c}\,\frac{dm(\z)}{\rho^q(\z)}.
\end{equation}
We proceed with another application of \eqref{lcr} to obtain
\begin{equation}\label{e140}
\begin{split}
\int_{E_t^c}\,\frac{dm(\z)}{\rho^q(\z)} = & q\int_0^\infty
y^{q-1}\left|\left\{\z\in\T:\,t<\rho(\z)<1/y\right\}\right|\,dy
\\ = & 2^{-q}- \frac{\nu_E(t)}{t^q}+q\int_t^2\,\frac{\nu_E(s)}{s^{q+1}}\,ds\,.
\end{split}
\end{equation}
Hence, the following bound holds
\begin{equation}\label{e15}
\int_{\O} G_{\O}(0,\l)\,d\mu(\l)\le 6\cdot
2^{-q}+6q\int_t^2\,\frac{\nu_E(s)}{s^{q+1}}\,ds\,.
\end{equation}

Our next step is to obtain a lower bound for the Green function
$G_\O(0,\cdot)$ in a smaller domain $\O_{\tau}$ with
$\tau=(6\pi+3)t$. We assume now that $t\in (0, \frac{1}{6\pi+3})$,
and so $\tau\in (0,1)$. For $z\in\G$ there is $\z'\in E$ such that
$|z-\z'|=t$, so $t\ge 1-|z|$, and hence
$$ \log\frac1{|z|}\le\log\frac1{1-t}\le \frac{3t}2\,.$$
The harmonic function $h(0,\cdot)$ from \eqref{G} does not exceed
$3t/2$ on $\G$ and equals zero on $E_t^c$, so
\begin{equation}\label{e16}
h(0,\l)\le\frac{3t}2\,, \qquad \l\in\O. \end{equation}

We will distinguish two situations for $\l$.

\noindent Let $|\l|\le 1-2t$. Then by \eqref{e16}
 $$
G_\O(0,\l)=\log\frac1{|\l|}-h(0,\l)\ge 1-|\l|-\frac{3t}2\ge
\frac{1-|\l|}4.
 $$
 Let $|\l|>1-2t$. This is where the restriction
$\l\in\O_\tau$ is essential. Let $\o(\l,E_t,\D)$ be the harmonic measure of the set
$E_t$. We put $g(\l)=9t\o(\l,E_t,\D)$. The same argument as above implies
$\o(\l,E_t,\D)\ge 1/6$ for $\l\in\G$, so
$$ h(0,\l)\le g(\l), \quad \l\in \partial\O \ \ \Rightarrow
\ \ \ h(0,\l)\le g(\l), \quad \l\in\O, $$ and
 $$
G_\O(0,\l)\ge 1-|\l|-g(\l).
 $$
We will find the upper bound for $g(\l)$ in $\O_\tau$. For
$\l=|\l|e^{i\th}$ write
$$
g(\l)=9t\int_{E_t}\,\frac{1-|\l|^2}{|\z-\l|^2}\,dm(\z)=\frac{9t(1-|\l|^2)}{2\pi}
\int_{e^{i\p}\in E_t}\,\frac{d\p}{(1-|\l|)^2+4|\l|\,\sin^2\frac{\p-\th}2}\,.
$$
Since $\rho(\l)\ge \tau$ and $1-|\l|<2t$, we have for some $\z_1\in
E$
 \begin{eqnarray*}
\rho(e^{i\th}) &=& |e^{i\th}-\z_1|=|\l-\z_1+e^{i\th}-\l| \\ &\ge&
(6\pi+3)t-(1-|\l|)>(6\pi+1)t.
 \end{eqnarray*}
 Next, for
$e^{i\p}\in E_t$ there is $\z_2\in E$ with
$|e^{i\p}-\z_2|=\rho(e^{i\p})\le t$, so
 $$
|\th-\p|\ge2\left|\sin\frac{\th-\p}2\right| = |e^{i\th}-\z_2+\z_2-e^{i\p}|\ge
\rho(e^{i\th})-t \ge 6\pi t.
 $$
Therefore, as $|\l|>1-2t>9/10$, we come to
 $$
g(\l) \le \frac{5\pi t(1-|\l|)}{2}\,\int_{6\pi t\le |\p-\th|\le\pi}\,
\frac{d\p}{(\p-\th)^2}<\frac{5(1-|\l|)}6.
 $$
Finally, the following lower bound holds for all $\l\in\O_\tau$
$$ G_\O(0,\l)\ge\frac{1-|\l|}6\,, $$
and we can continue \eqref{e15} as
\begin{equation}\label{e19}
\int_{\O_\tau} (1-|\l|)\,d\mu(\l)\le 36 \lp
2^{-q}+q\int_{\tau/(6\pi+3)}^2\,\frac{\nu_E(s)}{s^{q+1}}\,ds\rp.
\end{equation}
The latter holds for each $\tau\in (0,1)$.

Define $\ti\rho:=\min\{\rho,1\}$. It is easily checked that
\begin{equation}\label{e19.1}
\ti\rho(\l)\le\rho(\l)\le 2\ti\rho(\l) \quad \forall \l\in\D,
\end{equation}
and $\{\l:\rho(\l)>\tau\}=\{\l:\ti\rho(\l)>\tau\}$ for $0<\tau<1$.
\eqref{lcr} applied to the measure $d\sigma=(1-|\l|)d\mu$ and the
function $\ti\rho\le 1$ gives
$$
\int_\D\,(\ti\rho(\l))^{q-\a}d\s(\l)=(q-\a)\int_0^1
\tau^{q-\a-1}\,d\tau\,\int_{\{\rho>\tau\}}\,d\s(\l). $$ Since
$\rho(z)\le 2\rho(\tau z)$ for $z\in\D$ and $0<\tau<1$ then
$\rho(\tau z)>\tau/2$, as long as $\rho(z)>\tau$, so the whole
interval $[0,z]$ belongs to the set $\{\rho>\tau/2\}$. The latter
means that $\{\rho>\tau\}\subset\O_{\tau/2}$. Hence by using
\eqref{e19.1} we have
$$
\int_\D\,(\rho(\l))^{q-\a}d\s(\l)\le 2^{q-\a}(q-\a)\int_0^1
\tau^{q-\a-1}\,d\tau\,\int_{\O_{\tau/2}}\,(1-|\l|)\,d\mu(\l), $$ and
\eqref{e19} leads to \eqref{blas4}, as claimed.

\medskip

To prove \eqref{blas1}, \eqref{blas3} in the case $v(0)>-\infty$,
one can apply \eqref{blas2}, (respectively, \eqref{blas4}) to the
function $v_1(z)=(v(z)-v(0))(1+2^q|v(0)|)^{-1}$. Indeed,
$v_1(z)\le\rho^{-q}(z)$
 for all $z\in\bd$, and  the Riesz measure of $v_1$ coincides with the Riesz measure of $v$
up to a constant factor (which depends on $v(0)$).

If $v(0)=-\infty$, take the harmonic function $h$ in the disk
$\{|z|<1/2\}$  such that $v=h$ for $|z|=1/2$, and put
 $$
v_1(z)=\left\{\begin{array}{ccc}\max(v(z),h(z))&\hbox{for}&|z|<1/2,
\\ v(z)&\hbox{for}& |z|\ge1/2.\end{array}\right.
 $$
Clearly, $v_1(z)\le v(z)\le\rho^{-q}(z,E)$. Moreover, $v_1$ is
subharmonic in $\bd$ (see, e.g., \cite[Theorem 2.4.5]{ran}), and the
restriction of its Riesz measure $\mu_1$ on the set
$\{z\in\bd:\,|z|>1/2\}$ is the same as one for $\mu$. Therefore, the
integrals in \eqref{blas1}, \eqref{blas3} taken for $\mu$ and
$\mu_1$ differ by a bounded term. Since $v_1(0)\neq\infty$, the
result follows. The proof is complete. \bs

\medskip

{\it Proof of Theorem \ref{t1.1}}. As it was mentioned above, the
function $v_0=\rho^{-q}$ is subharmonic in $\bd$. We invoke again
the harmonic function $V$ \eqref{harm}.  It is clear that $V=v_0$ on
$E_t^c$, and $V\le t^{-q}=v_0$ on $\G$. Put $u_0$ for the least
harmonic majorant of $v_0$ on $\O=\O_t$, which is defined above. We
have
$$ V(\z)\le v_0(\z)\le u_0(\z), \quad \z\in\partial\O \ \ \ \Rightarrow
\ \ \ V(z)\le u_0(z), \quad z\in\O. $$
 Since $h(0,\cdot)\ge0$ and $v_0(0)=1$, then by
\eqref{Rie} for $z=0$
\begin{equation}\label{e110}
V(0)\le u_0(0)\le \int_\O\,\log\frac1{|\l|}\,d\mu_0(\l)+1.
\end{equation}
For $V(0)$ we see (cf. \eqref{e14}--\eqref{e140}) that
\begin{equation}\label{e111}
 V(0)=2^{-q}+q\int_t^2\,
\frac{\nu_E(s)}{s^{q+1}}\,ds\,. \end{equation}
 Thereby, in the case $q=\a$ we have
 $$
\int_\D\,\log\frac1{|\l|}\,d\mu_0(\l)=+\infty.
 $$
Let now $q>\a$. We apply \eqref{lcr} with
$d\sigma=\log|\l|^{-1}\,d\mu_0$, keeping in mind
\eqref{e110}--\eqref{e111}
\begin{eqnarray*}
& &\int_\D \,\bigl(\rho(\l,E)\bigr)^{q-\a}d\sigma(\l) =
(q-\a)\int_0^2\,y^{q-\a-1}\sigma(\{\l:\rho>y\})\,dy \\
&\ge&
(q-\a)\int_0^1\,y^{q-\a-1}\,dy\,\int_{\O_y}\,\log\frac1{|\l|}\,d\mu_0(\l)\\
&\ge&
(q-\a)\int_0^1\,y^{q-\a-1}\left\{2^{-q}-1+q\int_y^1\,\frac{\nu_E(s)}{s^{q+1}}
\,ds\right\}\,dy \\ &=&
q\int_0^1\,\frac{\nu_E(s)}{s^{\a+1}}\,ds+2^{-q}-1,
\end{eqnarray*}
and hence
$$
\int_\D\log\frac1{|\l|}\bigl(\rho(\l,E)\bigr)^{q-\a}~d\mu_0(\l)=+\infty.
$$
\smallskip

To go over to \eqref{c} we note that in view of \eqref{Rie} and
\eqref{G} the function $\rho^{-q}$ in the disk $\{z:|z|<1/2\}$
differs from $\int_{|\l|<1/2}\log(1/|\l-z|)~d\mu_0(\l)$ by a
harmonic function. Hence, the latter integral is finite at $z=0$,
and \eqref{c} follows. The proof is complete. \bs

\smallskip

 The following example shows that the result in Theorem
\ref{t2} is optimal as well.

{\bf Example}. Let $z_n=1-1/(n+1)$, \ $n=1,2,\ldots$.  Consider an infinite canonical
product (see, e.g., \cite{hkz}, p. 132)
 $$
 f(z):=\prod_{n=1}^\infty P\lp\frac{1-z_n^2}{1-z_nz}\rp,\qquad
 P(z)=(1-z)\,e^z.
  $$
It follows from the Taylor expansion of $P$ that $|1-P(z)|<4|z|^2$
for $|z|\le2$. Since $|1-z_n^2|<2|1-z_nz|$ for $z\in\bd$, we get
 \begin{equation}\label{b}
\log|f(z)|\le \sum_{n=1}^\infty \left|1-P\lp\frac{1-z_n^2}{1-z_nz}\rp\right|
 <4\sum_{n=1}^\infty \left|\frac{1-z_n^2}{1-z_nz}\right|^2.
 \end{equation}
 Furthermore, for $z\in\bd$
  $$
\left|\frac{1-z_n^2}{1-z_nz}\right|\le \frac{2}{|n(1-z)+1|}\le
\frac{2\,\sqrt{2}}{n|1-z|+1}
  $$
Hence for $n(z)=[|1-z|^{-1}]+1\in\N$  the right hand side of
(\ref{b}) does not exceed
 $$
 32\,n(z)+\frac{32}{|1-z|^2}\sum_{n=n(z)+1}^\infty\frac1{n^2}\le \frac{64}{|1-z|}
 +\frac{32}{n(z)|1-z|^2}\le
\frac{96}{|1-z|}\,.
 $$
 So $f$ satisfies \eqref{e010} with $E=\{1\}$, $q=1$, and the choice
 $\a=1-\ep$ is optimal.

 The growth of certain canonical products in the unit disk was
 studied in \cite{G, G1}.

\section{Discrete spectrum of contraction operators}\label{s2}

The relation between the spectral theory of non-selfadjoint operators and the theory
of analytic functions via perturbation determinants is well established due to the
works of I.~Gohberg, M.~Krein, V.~Matsaev and others. The variants of this approach
were used recently in \cite{dk, dhk}.

The results from Section 2, in particular Theorem \ref{t2}, can be applied to the
operator theory in the same fashion as it is done in \cite{bgk, dk, dhk}.

Let $U$ be a unitary operator in the Hilbert space $\ch$, $\s(U)=E$
its spectrum. The resolvent $R_U(z)=(U-z)^{-1}$ is an analytic in
$\bd$ operator-function, and (see, e.g., \cite[Section V.3.8]{kato})
$$ \|R_U(z)\|=\rho^{-1}(z,\s(U)) \qquad z\in\bd. $$

A bounded, linear operator $T$ in $\ch$ is called a {\it
contraction} if $\|T\|\le1$, $\|\cdot\|$ is the standard operator
norm. Assume that $T-U$ is a compact operator. If $E\not=\bt$, then
by a general result from the perturbation theory
$\s(T)=E\cup\s_d(T)$, where the discrete component of the spectrum
$\s_d(T)\subset\ovl\bd\backslash E$ is at most countable set with
all possible accumulation points in $E$, and each $z_k\in\s_d(T)$ is
an eigenvalue of finite algebraic multiplicity.

Our goal here is to obtain some quantitative bound for the rate of convergence of
$z_k$ to $E$ in the case, when $T-U\in\css_q$, $q>0$, the Schatten-von Neumann classes
of compact operators, i.e., the sequence of eigenvalues of the operator $|T-U|$
belongs to $\ell^q$. The norm in $\css_q$ is denoted by $\|\cdot\|_q$, and the
standard operator norm $\|\cdot\|=\|\cdot\|_\infty$. The family $\{\css_q\}$ is known
to be nested: $\css_{q_1}\subset\css_{q_2}$ for $0<q_1<q_2\le\infty$, and
$\|A\|_{q_2}\le\|A\|_{q_1}$. If $A\in\css_q$, $q\ge1$, and $S$ is a bounded linear
operator then $AS\in\css_q$ and $\|AS\|_q\le \|A\|_q\,\|S\|$. An extensive information
on the subject is Gohberg--Krein \cite{gk}, Simon \cite{si1}, and Dunford--Schwartz
\cite[Chapter XI.9]{dsh}.

The key analytic tool is the so called renormalized determinant
$\det_n(I+A)$, defined for $A\in\css_n$ and $n=1,2,\ldots$. In
particular, we need the following continuity property (cf.
\cite[Theorem 9.2, (c)]{si1})
\begin{equation}\label{rdet}
|\det{}_n(I+A)-\det{}_n(I+B)|\le \|A-B\|_n\exp\lp
C_n(1+\|A\|_n^n+\|B\|_n^n)\rp, \end{equation} a constant $C_n$
depends only on $n$. If $T-U\in\css_n$, then the perturbation
determinant
$$ u_n(z):=\det{}_n(T-z)(U-z)^{-1}=\det{}_n\lp I+(T-U)(U-z)^{-1}\rp $$
is a well defined and analytic function in $\bd$, and its zero set
in $\D$ is $Z_{u_n}=\s_d(T)\cap\bd$, with the order of each zero
equal the algebraic multiplicity of the eigenvalue.

\begin{theorem}\label{t210}
Let $U$ be a unitary operator in $\ch$,  $\s(U)=E\ne\T$, and
$I(\a,E)<\infty$ in \eqref{aset} for some $\a\in\R$. Let a
contraction $T$ be a $\css_q$-perturbation of $U$: $T-U\in\css_q$,
$q>0$. Then
\begin{equation}\label{e2.1} \sum_{z\in \s_d(T)}
(1-|z|)\rho(z,E)^{(q-\a)_+}<\infty. \end{equation}
 Moreover, there is a constant $\d=\d_q$, $0<\d_q<1$, which depends only on
$q$, such that for $\|T-U\|_q<\d_q$ the inequality holds
\begin{equation}\label{e2.0}
\sum_{z\in \s_d(T)} (1-|z|)\rho(z,E)^{(q-\a)_+}\le
C(q,E)\,\max(\|T-U\|_q,\|T-U\|_q^q).
\end{equation}
\end{theorem}

{\sl Remark}. Obviously, the eigenvalues of $T$ with $|z|=1$ do not enter the left
hand side, so we claim nothing about this part of the discrete spectrum. However, this
part is missing as long as $T$ is a completely non-unitary contraction.

\smallskip

{\it Proof}. Take a positive integer $n$ from $q\le n<q+1$, so
$T-U\in\css_n$. By \cite[Lemma XI.9.22,(d)]{dsh}
\begin{eqnarray*}
|u_n(z)| &\le& \exp\lp c_q\|(T-U)(U-z)^{-1}\|_q^q\rp \\
&\le& \exp\lp c_q\|T-U\|_q^q\,\|(U-z)^{-1}\|^q\rp \\
&\le& \exp\lp c_q\|T-U\|_q^q\rho^{-q}(z,E)\rp.
\end{eqnarray*}
An application of the first statement of Theorem \ref{t2} proves
\eqref{e2.1}.

Throughout the rest of the proof $c_q$ stands for different positive
constants which depend only on $q$. Since $\|T-U\|<1$, then $T$ is
invertible, and
 $$
 u_n(0)=\det{}_n\lp I+(T-U)U^{-1}\rp\not=0.
 $$
Put $f(z)=u_n(z)/u_n(0)$, $|f(0)|=1$, so that
 $$
 \log |f(z)|\le \frac{c_q\|T-U\|_q^q}{\rho^q(z,E)}
+\log\frac1{|u_n(0)|}.
 $$
Bound \eqref{rdet} with $A=(T-U)U^{-1}$, $B=0$ gives in view of the
monotonicity of norms
 $$
|u_n(0)-1|\le \|T-U\|_q\exp\lp C'_q(1+\|T-U\|_q^n)\rp<\frac12,
 $$
 as soon as $\d_q$ is small enough. Hence
 $$
\log\frac1{|u_n(0)|} \le \log 4\,|u_n(0)-1|\le c_q\,\|T-U\|_q.
 $$
 An application of the second statement of Theorem \ref{t2} gives
 \eqref{e2.0}. \bs

\end{document}